\documentclass[10pt]{amsart}
\theoremstyle{plain}
\usepackage{amsmath, amssymb, amscd,graphicx}
\newtheorem{lem}{Lemma}[section]
\newtheorem{prop}[lem]{Proposition}
\newtheorem{thrm}[lem]{Theorem}
\newtheorem{cor}[lem]{Corollary}
\newtheorem{dfn}[lem] {Definition}

\newtheorem{tm}{Theorem}
\newtheorem*{conj}{Conjecture}
\newtheorem{crrr}{Corollary}

\newcommand{\Szabo}{{Szab{\'o} }}
\newcommand{\smax}{{s_{\text{max}}}}
\newcommand{\smin}{{s_{\text{min}}}}
\newcommand{\co}{{\colon}}
\newcommand{\s}{{\mathfrak{s}}}
\renewcommand{\b}{{\bf b}}
\renewcommand{\a}{{\bf a}}
\newcommand{\q}{{q}}
\newcommand{\x}{{\bf v}}
\newcommand{\ckh}{{CKh}}
\newcommand{\kh}{{Kh}}

\newcommand{\R}{{\bf{R}}}

\newcommand{\Z}{{\bf{Z}}}
\newcommand{\Q}{{\bf{Q}}}

\newcommand{\into}{{\hookrightarrow}}

\newcommand{\gr}{{{\mathrm{gr} }}}

\newcommand{\A}{{{\mathcal A }}}

\newcommand{\ts}{{{\thinspace}}}

\title{Khovanov homology and the slice genus}
\author{Jacob Rasmussen}
\address{Princeton University Dept. of Mathematics, Princeton, NJ 08544}
\email{jrasmus@math.princeton.edu}
\thanks{The author is supported by an NSF Postdoctoral fellowship.}

\begin{document}

\begin{abstract}
We use Lee's work on the Khovanov homology to
define a knot invariant \(s\). We show that \(s(K)\) is a concordance
invariant and that it provides a lower
bound for the slice genus of \(K\). As a corollary, we give a purely
combinatorial proof of the Milnor conjecture. 
\end{abstract}

\maketitle

\section{Introduction}

In \cite{Khovanov}, Khovanov introduced an invariant of knots and
links, now widely known as the Khovanov homology. This invariant takes
the form of a graded homology theory \( \kh (L) \), whose graded Euler
characteristic is the unnormalized Jones polynomial of \(L\). In
\cite{ESL2}, Lee showed that \( \kh (L) \) is naturally viewed as the
\(E_2\) term of a spectral sequence which converges to \( \Q \oplus
\Q\). In this paper, we use this spectral sequence to define a knot
invariant \(s(K) \). The definition of \(s(K) \) was motivated by a
similar invariant \( \tau (K) \) which is defined using knot Floer
homology \cite{OS10}, \cite{thesis}. In fact, the similarities between
the two invariants extend far beyond their manner of definition.

Our main result is that the invariant \(s\) gives a
lower bound for the slice genus:
\begin{tm}
\label{Thm:One}
\begin{equation*}
|s(K)| \leq 2 g_*(K)
\end{equation*}
where \( g_*(K) \) denotes the slice genus. 
\end{tm}
\noindent In fact, 
\begin{tm}
\label{Thm:Two}
The map \(s\) induces a homomorphism from \(\text{Conc}(S^3) \) to \(
\Z\), where \( \text{Conc}(S^3) \) denotes the concordance group of
knots in \(S^3\). 
\end{tm}

\noindent For alternating knots, \(s(K)\) does not provide any new information
about \(g_*(K)\):
\begin{tm}
\label{Thm:Three}
If \(K\) is an alternating knot, then \(s(K) \) is equal to the
classical knot signature \( \sigma (K) \). 
\end{tm}
There is, however, a class of knots for which \(s(K) \) gives much
better --- indeed, sharp --- information. We say that a knot is
{\it positive} if it admits a planar diagram with all positive
crossings. 
\begin{tm}
\label{Thm:Four}
If \(K\) is a positive knot, 
\begin{equation*}
s(K) = 2 g_*(K) = 2 g(K)
\end{equation*}
where \(g(K) \) is the ordinary genus of \(K\). 
\end{tm}

As a corollary, we get a Khovanov homology proof of following result,
 which was first proved by Kronheimer and Mrowka using gauge theory
\cite{KMMilnor}:

\begin{crrr}
(The Milnor Conjecture) The slice genus of the \((p,q)\) torus knot is 
\((p-1)(q-1)/2\). 
\end{crrr}

As the reader familiar with knot Floer homology will have already
noted, the theorems above all hold with  \(2 \tau (K) \) in place of
  \(s (K) \). (See \cite{OS10} for the first three, and
  \cite{Livingston} and \cite{Rudolph} for
 the final one.) Indeed, the equality \( s(K) = 2 \tau (K) \) holds
in all cases for which the author knows the value of \( \tau (K) \). 
Based on these observations, we make the following (perhaps
optimistic) 
\begin{conj}
For any knot \(K \subset S^3\), \(s(K) = 2 \tau (K).\)
\end{conj}

Readers familiar with the Khovanov homology may also have observed that the
notation \(s(K) \) has already been used by Bar-Natan \cite{DBN} to
describe an apparent knot invariant which appears in one of his
``phenomenological conjectures.'' This is no coincidence. Indeed, the
author's interest in the subject was first aroused by the observation
that Bar-Natan's \(s\) appeared to  give a lower bound for the slice
genus. Although we are unable to prove that the \(s(K)\) defined here
is the same as that determined by Bar-Natan's conjecture,
we do give a fairly general condition (at least for small knots) 
under which the two agree. 

The remainder of the paper is organized as follows. In section 2, we
review the Khovanov complex and Lee's construction of a spectral
sequence from it. In section 3, we define \(s\) and show that it
behaves nicely with respect to the structure of the concordance
group. Section 4 is devoted to the proof of Theorem~\ref{Thm:One}. 
In section 5, we prove Theorems~\ref{Thm:Three} and \ref{Thm:Four},
and discuss the relationship between \( s(K) \) and \( \tau (K) \) in
more detail. Finally, section 6 contains proofs of some technical
results establishing the invariance of Lee's spectral sequence, which are
are needed in section 2. 

Finally,  we take this opportunity to fix two conventions which we will
use throughout. First, we will always work with \( \Q\) coefficients. 
 Although Khovanov's complex can be defined with coefficients in
\( \Z \),  Lee's theorem (Theorem~\ref{Thm:Lee2}) does not hold in
this context. Second, we will often abuse our notation, letting \(L\)
refer both to a planar diagram of a link and to the underlying link
itself. The reader should have little trouble determining from context
which meaning is intended. 

\noindent{\bf Acknowledgements:} The author would like to thank Peter
Kronheimer, Peter Ozsv{\'a}th, and Zoltan \Szabo for many helpful
conversations and for encouraging him to pursue this problem. 

\section{Review of Khovanov homology}
\label{Sec:Review}
In this section, we briefly recall the construction of the Khovanov
complex  \cite{Khovanov} and Lee's extension of it
\cite{ESL2}.

\subsection{The cube of resolutions}
\begin{figure}
\includegraphics{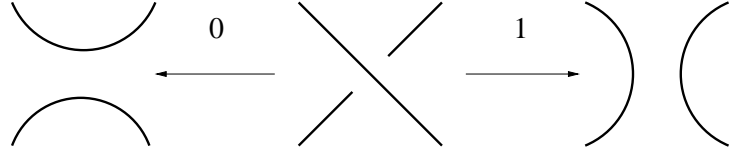}
\caption{\label{Fig:Resolutions} \(0\)-- and \(1\)--resolutions of
  a crossing.}
\end{figure}

Given a link diagram \(L\) with crossings labeled \(1\)
through \(k\), we can form the cube of all
possible resolutions of \(L\). This is a \(k\)-dimensional cube with its
vertices and edges decorated by 1-manifolds and
cobordisms between them. More specifically, each crossing of
\(L\) can be resolved in two different ways,
as illustrated in Figure~\ref{Fig:Resolutions}. To each vertex \(v\)
of the cube \( [0,1]^k\), we associate the planar diagram \(D_v\)
obtained by resolving the \(i\)-th crossing of \(L\) according to
 the \(i\)-th coordinate of \(v\). Then \(D_v\) is a planar
diagram without crossings, so it is a disjoint union of circles. 

Let \(e\) be an edge of the cube. The coordinates of its two ends
differ in one component --- say the \(l\)-th. 
 We call the end
which has a \(0\) in this component the {\it initial end}, and denote
it by \(v_e(0)\). The other end is called the {\it terminal end},
written \(v_e(1)\). We assign to
\(e\) the cobordism \(S_e: D_{v_e(0)} \to D_{v_e(1)} \), which is a
product cobordism except in a neighborhood of the \(l\)-th crossing,
where it is the obvious saddle cobordism between the \(0\) and
\(1\)-resolutions.

The Khovanov complex is constructed  by applying a  \(1+1\)
dimensional TQFT  \(\A\) to the cube of resolutions. 
In other words, one replaces each vertex \(v\) by the group \(\A
(D_v)\), and each edge \(e\) by the map \( \A (S_e) \). The underlying 
group of \(\ckh(L)\) is the direct sum of the groups \(\A (D_v)\) for all
vertices \(v\), and the differential on the summand \(\A (D_v) \) is a sum
of the maps \(\A(S_e) \) for all edges \(e\) which have  \(v\) as
their initial end. More precisely, for \(x \in \A (D_v) \)
\begin{equation}
\label{Eq:Diff}
d(x) = \sum_{i=1}^{c_0(v)} (-1)^{s(e_i)}\A(S_{e_i}).
\end{equation}
Here \(c_0(v) \) is the number of crossings in \(v\) which have a
  \(0\)--resolution, and \(e_i\) is the edge which  corresponds to
  changing the \(i\)-th such crossing to a \(1\)-resolution. 
The signs \((-1)^{s(e_i)}\) are chosen in such a way that
  \(d^2=0\). (There are many different ways to do this, but it is easy
  to see that they all give rise to isomorphic chain complexes.) 
The {\it homological grading} of an element \(x \in \A (D_v) \) is
  defined to be
\( \gr (v) = |v| - n_-\), where \(|v|\) is the number of \(1\)'s in
  the coordinates of \(v\) and \(n_-\) is the number of negative
  crossings in the diagram for \(L\). 
Note that  \(d\) increases the homological grading
  by \(1\) --- strictly speaking, the Khovanov homology is a
  cohomology theory!

\subsection{Khovanov's TQFT}
We now give a more explicit description of the TQFT \(\A\). 
 Let \(V \) be a vector space spanned by two elements,
\(\x_+\) and \(\x_-\). The vector space associated by \(\A\) to a
single circle is defined to be \(V\), so that if \(D\) is a diagram
composed of \(n\) disjoint circles, \(\A(D) = V^{\otimes n} \). Thus
we can think of \(\ckh (L) \) as being the vector space spanned by the
space of ``states'' for \(L\), where a state consists of a complete
resolution of \(L\),  together with a labeling of each component
of the resolution by either \(\x_+\) or \(\x_-\).
 
 The cobordisms \(S_e\) come in two
forms: either two circles can merge into one, or one can split into
two. In the first case, \(\A(S_e) \) is given by a map \(m \co
V^{\otimes 2} \to V\), where the two factors in the tensor product
correspond to the labels on the two circles that merge, and the copy
of \(V\) in the image corresponds to the label on the single resulting
circle. Likewise, in the second case, \( \A(S_e) \) is given by a map 
\(\Delta \co V \to V^{\otimes 2}\). The formulas for these maps are
\begin{align}
\label{Eq:Cobordisms1}
m(\x_+ \otimes \x_+) & = \x_+   \qquad &\Delta (\x_+ ) & = \x_+
\otimes \x_- + \x_- \otimes \x_+ \notag \\
m(\x_+\otimes \x_-) &= m(\x_- \otimes \x_+ ) = \x_- \qquad & 
\Delta  (\x_-) & = \x_- \otimes \x_- \\
m(\x_- \otimes \x_-) & = 0. & & \notag
\end{align}

 For reference, we also record two other
 maps \( \iota \) and \(\epsilon\) used to define \( \A\).
 These maps are not needed at the moment, but they make an
 appearance in section~\ref{Sec:Cobordisms}
 when we study cobordisms. 
Corresponding to the addition of a \(0\)-handle (the birth of a circle
 in a diagram), there is a map \( \iota \co \Q \to V \), and
 corresponding to the addition of a two handle (the death of a circle)
there is a map \( \epsilon \co V \to \Q \). These maps are given by 
\begin{align*}
\epsilon (\x_-) & = 1 \qquad & \iota (1) & = \x_+ \\
\epsilon (\x_+) & = 0. & &
\end{align*}

 \(\A\) is especially nice because it
 is a graded TQFT. We define a grading \(p\) on \(V\) by setting 
\(p(\x_{\pm}) = \pm 1\) and extend it to \(V^{\otimes n} \) by 
\(p(v_1\otimes v_2 \otimes \ldots \otimes v_n) = p (v_1) + p(v_2) +
\ldots  p (v_n) \). Then it is easy to see that if \({\bf v}\) is a homogenous
element of \(V^{\otimes n} \), \(p(S_e({\bf v})) = p({\bf v} ) - 1. \)
Next, we define a grading \(\q\) on \( \ckh (L) \) by 
\(\q ({\bf v}) = p({\bf v}) + \gr ({\bf v}) + n_+ -n_- \), where \(n_\pm\)
 are the number of positive and negative crossings in the diagram
 \(L\). (The term \(n_+ - n_-\) is included so that the \(q\)-grading
 remains invariant for different diagrams of the same knot.)
Then \(\q (d({\bf v})) = \q({\bf v}) \), so
\(\ckh (L) \) splits into a direct sum of complexes,
one for each \(q\) grading. In fact, its graded Euler characteristic
is the unnormalized Jones polynomial of \(L\), but we will not make
 use of this here. 

In \cite{Khovanov}, Khovanov proves that the homology of \(\ckh(L) \)
(thought of as a bigraded group) is an invariant of the underlying
link \(L\). We denote this homology group by \(\kh (L)\). 

\subsection{Lee's TQFT}
In \cite{ESL2}, Lee considers a similar construction, but with 
another TQFT \(\A'\) in place of \(\A\). The underlying vector spaces
for these two TQFT's are the same, but the maps \( m' \co V \otimes 
V \to V \) and 
\( \Delta' \co V \to V \otimes V \) induced by cobordisms 
are slightly different. They are given by 
\begin{align}
\label{Eq:Cobordisms2}
m'(\x_+ \otimes \x_+) & = m'(\x_- \otimes \x_-)= \x_+   
\qquad &\Delta ' (\x_+ ) & = \x_+ \otimes \x_- + \x_- \otimes \x_+ \\
m'(\x_+\otimes \x_-) &= m'(\x_- \otimes \x_+ ) = \x_- \qquad & 
\Delta ' (\x_-) & = \x_- \otimes \x_- +\x_+ \otimes \x_+. \notag
\end{align}
(The maps \( \iota \) and \( \epsilon \) corresponding to the addition
of \(0\) and \(2\)-handles are the same as for \( \A \).)
We denote the resulting complex by \(\ckh'(L)\) and its homology by 
\( \kh '(L) \). 

Using the obvious identification between the underlying groups of 
\(\ckh (L) \) and \(\ckh'(L) \), we can define a \(q\)-grading on the
latter group as well.
It is clear from equation~\ref{Eq:Cobordisms2} that this grading does
not behave quite so well with respect to the differential
\(d'\). Indeed, \(\Delta '(\x_-)\) is not even homogenous. It is easy
to see, however, that if \(\x  \in \ckh ' (L) \) is a homogenous
element, then the \(q\)-grading of every monomial in \( d' (\x)
\) is greater than or equal to the \(q\)-grading of \(\x\). In other
words, the \(q\)-grading defines a filtration on the complex \(\ckh
'(L) \). This fact leads to the following theorem, which is implicit
in \cite{ESL2}:
\begin{thrm}
\label{Thm:Lee1}
There is a spectral sequence with \(E_2\) term  \(\kh(L) \) 
which converges to \( \kh ' (L) \). The \(E_2\) and higher terms of
this spectral sequence are invariants of the link \(L\). 
\end{thrm}

The first part of the theorem is more or less immediate from the
observations above. The filtration on \( \ckh' \) gives rise to a
spectral sequence converging to \( \kh'\). The differential on its
\(E_1\) term is the part of \( d'\) which preserves (rather than
raises) the \(q\)-grading. Comparing equations~\ref{Eq:Cobordisms1} and 
\ref{Eq:Cobordisms2}, we see that the \(E_1\) term is  the
complex \( \ckh \). 

To prove the second statement, we check that the spectral sequence is
invariant under the Reidemeister moves. Suppose \(L\) and
\(\tilde{L}\) are two diagrams related by the \(i\)-th Reidemeister
move. In \cite{ESL2},  Lee defines maps \( \rho_i ' \co \ckh ' (L) \to \ckh
' (\tilde{L}) \) which induce isomorphisms on homology. 
In section~\ref{Sec:Reidemeister}, we show that these maps induce
isomorphisms on  \(E_2\) terms of spectral
sequences, thus completing the proof of the theorem. 

\subsection{Calculation of \(\kh'\).}
 Lee's second major result is that the homology group
 \(\kh'(L)\) is surprisingly simple. To show this, she
introduces a new basis \(\{ \a,\b \} \) for  \(V\),
where \(\a = \x_- + \x_+ \) and \(\b = \x_- - \x_+ \). With respect to
this new basis, the maps \(m'\) and \(\Delta'\) are given by
\begin{align*}
m'(\a \otimes \a) & = 2\a   \qquad &\Delta ' (\a) & = \a \otimes \a \\
m'(\a \otimes \b) &= m'(\b \otimes \a) = 0 \qquad & 
\Delta ' (\b) & = \b \otimes \b \\
m'(\b \otimes \b) & = -2 \b & &
\end{align*}
and the maps \( \epsilon '\) and \(\iota '\) are given by 
\begin{align*}
\epsilon ' (\a) & = \epsilon' (\b) = 1 \qquad & \iota ' (1) & = (\a-\b)/2
\end{align*}
Using this basis, she proves

\begin{thrm}
\label{Thm:Lee2}
(Theorem 5.1 of \cite{ESL2})
\(\kh'(L) \) has rank \(2^n\), where \(n\) is the number of components
of \(L\). 
\end{thrm}

Indeed, Lee exhibits an explicit bijection between the set of possible
orientations for \(L\) and a set of generators of \( \kh'(L) \), which
we refer to as {\it canonical generators}. This bijection may be
described as follows. Given an orientation \(o\) of \(L\), let \(D_o\)
be the corresponding  oriented resolution. We label the circles in \(D_o\) with
\(\a\)'s and \(\b\)'s according to the following rule. To each circle \(C\)
we assign a mod 2 invariant, which is the mod 2 number of
circles in \(D_o\) which separate it from infinity. (In other words,
draw a ray in the plane from \(C\) to infinity, and take the number
of other times it intersects the other circles, mod 2.) To this
number, we add \(1\) if \(C\) has the standard (counterclockwise)
orientation, and \(0\) if it does not. Label \(C\) by \(\a\) if the resulting
invariant is \(0\), and by \(\b\) if it is \(1\). We denote the
resulting state by \( \s _o\).

The name ``canonical generator'' is justified by the following result,
whose proof is given in section~\ref{Sec:Reidemeister}. 
\begin{prop}
\label{Prop:ReidCan}
Suppose \(L\) and \(\tilde{L} \) are related by the \(i\)-th
 Reidemeister
move. Then an orientation \(o\) on \(L\) induces an orientation
\(\tilde{o}\) on \( \tilde{L}\), and \( \rho _{i*} ' ([\s_o])\) is
 a nonzero multiple of \( [\s_{\tilde{o}}] \). 
\end{prop}

\begin{figure}
\includegraphics{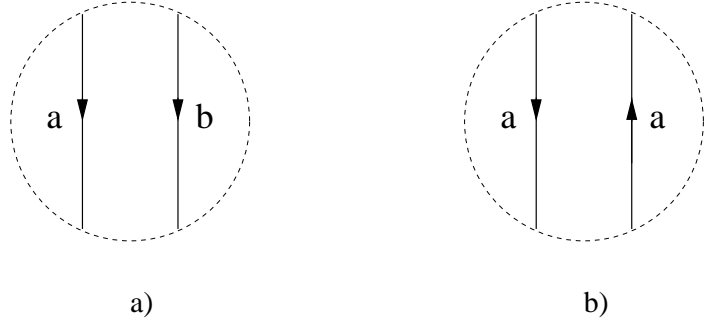}
\caption{\label{Fig:Segments} Local behavior of the state \(\s_o\).}
\end{figure}
 We end this section with an  elementary but important observation.
\begin{lem}
\label{Lem:CO}
(Coherent orientations) Suppose there is a region in the state diagram
    for \(\s_o\)
    containing exactly two segments, as shown in
    Figure~\ref{Fig:Segments}. Then either the orientations of the two are
    the same and the labels are different (like part {\it a} of the
    figure) or the orientations are different and the labels are the
    same (like part {\it b}). 
\end{lem} 

\begin{proof}
We consider three possible cases: either the two segments belong to
the same circle in \(D_o\), or they belong to two circles, one of
which is contained inside the other, or they belong to two circles,
neither of which is contained inside the other. In each case, it is
easy to verify that the claim holds. 
\end{proof}

\begin{cor}
If two circles in the state diagram for
\(\s_o\) share a crossing, they have different labels. 
\end{cor}

\section{Definition and Basic Properties of the Invariant}
\label{Sec:Def}
Let  \(K\) be  a knot in \(S^3\). By
Theorems~\ref{Thm:Lee1} and \ref{Thm:Lee2}, we know that there is a
spectral sequence associated to \(K\) which converges to \(\Q\oplus
\Q\). This spectral sequence is a relatively complicated object, but
we can  extract some simpler invariants of \(K\) from it.
 Let \(\smax \) and \(\smin\) (with \(\smax \geq \smin \))
be the \(q\)-gradings of the two surviving copies of \(\Q\) which
remain in the \(E_{\infty}\) term of the spectral sequence.
Like all \(q\)-gradings for a
knot, \( \smax \) and \( \smin \) are odd integers. Since the
isomorphism type of the spectral sequence is an invariant of \(K\),
\(\smax \) and \(\smin\) are invariants as well. 

Before making this definition formal, we digress  to
 establish some  terminology related to filtrations.
Suppose \(C\) is a chain complex. A {\it finite length filtration} of
 \(C\) is a sequence of subcomplexes
\begin{equation*}
 0 = C_{n} \subset C_{n-1} \subset C_{n-2} \subset \cdots \subset
C_m = C.
\end{equation*} 
To such a filtration, we associate a {\it grading} defined as
follows: \(x \in C\) has grading \(i\) if and only if 
\(x \in C_i\) but \( x \not \in C_{i-1} \). 
If \(f: C \to C' \) is a map between two filtered chain complexes, we
say that \(f\) {\it respects the filtration} if \(f(C_i) \subset C'_i
\). More generally, we say that \(f\) is a {\it filtered map of degree
 } \(k\) if \( f(C_i) \subset C_{i+k}'.\)

A  filtration \(\{C_i\}\)
 on \(C\) induces a filtration \( \{S_i\} \)  on \(H_*(C) \)
 defined as follows: a class \([x] \in H_*(C) \)
is in \(S_i\) if and only if has a representative which is an element
of \(C_i\). If \(f \co C \to C' \) is a filtered chain map of degree
 \(k\), then it is easy to see that the induced map 
\( f_* \co H_*(C) \to H_*(C') \) is also filtered of degree \(k\). 

A finite length filtration\(\{C_i\}\) on \(C\) induces a spectral
 sequence, which 
 converges to the associated graded group of the induced filtration 
\(\{S_i\} \). In other words, the group which survives at grading \(i\)
in the spectral sequence is naturally identified with
 the group \(S_i/S_{i+1}\). 

Let us denote by \(s\) the grading on \( \kh ' (K) \) induced by the
\(q\)-grading on \( \ckh ' (K) \). Then the informal definition above
is equivalent to 
\begin{dfn}
\begin{align*}
\smin (K) & = \min \{s(x) \ts | \ts x \in \kh' (K), x \neq 0 \} \\
\smax (K) & = \max \{s(x) \ts | \ts x \in \kh' (K), x \neq 0 \}
\end{align*}
\end{dfn}

\noindent Since \(\kh\) of the unknot \(U\)
 has rank two and is supported in
\(q\)-gradings \( \pm 1\), we have \( \smax (U) = 1\), \( \smin
(U) = -1\). 

Another proof that \( \smax \) and \( \smin \) are knot invariants
could be given using

\begin{prop}
\label{Prop:ReidFilt}
The maps \( \rho_{i*}'\) and \( (\rho_{i_*}')^{-1} \) both respect the
induced filtration \(s\) on \( \kh' \). 
\end{prop}
\noindent The proof may be found in section~\ref{Sec:Reidemeister}.

\subsection{The invariant \(s\)} Our first task in this section is to
prove

\begin{prop}
\label{Prop:MinMax}
\begin{equation*}
\smax(K) = \smin(K) + 2
\end{equation*}
\end{prop}

\noindent which justifies 

\begin{dfn}
\begin{equation*}
s(K) = \smax(K) -1 = \smin(K) + 1
\end{equation*}
\end{dfn}
\noindent Since \(\smax \) and \(\smin \) are odd, \( s (K) \) is always an
even integer. 

Before proving the proposition, we need some preliminary results.

\begin{lem}
Let \(n\) be the number of components of \(L\).
There is a direct sum decomposition \(\kh'(L) \cong \kh'_o(L) \oplus
\kh'_e(L) \), where \( \kh'_o (L) \) is generated by all states with
\( q\)-grading conguent to \(2+n \mod{4}\), and \( \kh'_e (L) \) is generated
by all states with \( q\)-grading congruent to \(n \mod{4}\).
 If \(o\) is
an orientation on \(L\),
 then \(\s_o + \s_{\overline{o}} \) is contained in one of the two
 summands, and \(  \s_o - \s_{\overline{o}} \) is contained in the other.
\end{lem}

\begin{proof}
Following Lee \cite{ESL2}, we write
\begin{align*}
m' & = m + \Phi_m \\
\Delta' & = \Delta + \Phi_\Delta
\end{align*}
where \( m \) and \( \Delta \) preserve the \(q\)-grading and \(
\Phi_m\) and \( \Phi_\Delta \) raise it by 4. This proves the first
statement. 

For the second, let \( \iota \co \ckh ' (L ) \to \ckh ' (L)
\) be the map which acts by the identity on \( \ckh' _e \) and by
multiplication by \( -1\) on \( \ckh' _o \). We claim that \( \iota
(\s_o) = \pm \s_{\overline{o}} \). To see this, we define
a new grading on \(V\) with respect to which \(\x_-\) has  grading 0 and
\(\x_+\) has grading 2. Let  \(i \co V \to V\) be given by \(i(\x_-) =
\x_-\), \(i(\x_+) = -\x_+\), so that  \(i (\a) = \b \) and \(i(\b) = \a
\). Then the induced map \(i ^{\otimes n} \co V^{\otimes n} \to V^{\otimes n}
\) acts as the identity on elements whose new grading is congruent
to \( 0 \mod{4} \) and as multiplication by \(-1\) on elements whose
new grading is congruent to \( 2 \mod{4} \). The new grading differs
from the \(q\)-grading on \(D_o\) by an overall shift, so 
\begin{equation*}
\iota (\s_o) = \pm i ^{\otimes n} (\s_o) = \pm \s_{\overline{o}}
\end{equation*}
It follows that \( \s_o + \iota (\s_o) = \s_o \pm \s_{\overline{o}}\) is
contained in one summand, while \( \s_o - \iota (\s_o) = 
\s_o \mp \s_{\overline{o}} \) is contained in the other. 
\end{proof}

\begin{cor}
\begin{equation*}
s(\s_o) = s(\s_{\overline{o}}) = \smin(K)
\end{equation*}
\end{cor}

\begin{cor}
\label{cor:Neq}
\( \smax (K) > \smin (K) \). 
\end{cor}

\begin{proof}
Since \( \ckh' (K) \) decomposes as a direct sum, its affiliated
spectral sequence decomposes too. The homology of each summand is
\(\Q\), so each must account for one of the surviving terms in the
spectral sequence. The two summands are supported in different
\(q\)-gradings, so the surviving terms must have different
\(q\)-gradings as well. 
\end{proof}


\begin{lem}
\label{Lem:SumSeq}
For knots \(K_1\), \(K_2\), there is a short exact sequence
\begin{equation*}
\begin{CD}
0 @>>> \kh' (K_1 \# K_2 ) @>{p_*}>> \kh'(K_1) \otimes \kh'(K_2) @>{\partial}>> 
\kh' (K_1 \# K_2 ) @>>> 0 
\end{CD}
\end{equation*}
The maps \(p_*\) and \(\partial\) are filtered of degree \(-1\). 
\end{lem}

\begin{figure}
\includegraphics{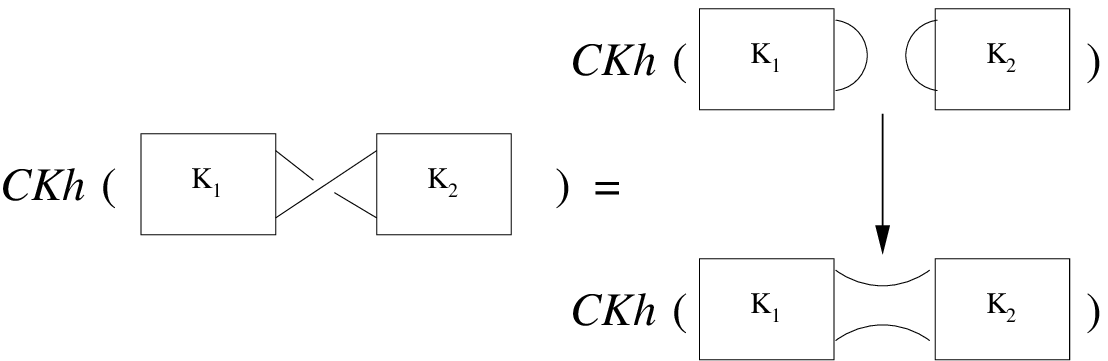}
\caption{\label{Fig:Sum} A short exact sequence for \( \ckh '(K_1 \#
  K_2) \). }
\end{figure}

\begin{proof}
Consider the diagram 
 for \(K_1 \# K_2\) shown in Figure~\ref{Fig:Sum}. From it, we get
 a short exact sequence
\begin{equation*}
\begin{CD}
0 @>>> \ckh' (D_1)\{1\} @>>> \ckh' (D_2) @>{p}>> \ckh' (D_3) @>>> 0
\end{CD}
\end{equation*}
where \(D_1\) and \(D_2\) are both diagrams for \(K_1 \# K_2 \),
and  \(D_3\) is a diagram for the disjoint union \(K_1 \coprod K_2 \).
Since \(\kh' (K_1 \# K_2) \) has rank two and 
\(\kh' (K_1 \coprod K_2) \cong \kh' (K_1) \otimes \kh' (K_2) \) has
rank four, the resulting long exact sequence must split, giving the
 short exact sequence of the lemma. It is clear that the maps \(p_*\) and \(
\partial\) are filtered of some degree, which can be worked out
by considering (for example) the case \(K_1 = K_2 = U\). 
\end{proof}

\begin{proof} (of Proposition~\ref{Prop:MinMax}.)
Consider the exact sequence of the previous lemma with \(K_1 = K\) and 
\(K_2\) the unknot. Denote the canonical generators of \(K\) by \(
\s_a\) and \( \s_b\), according to their label near the connected sum
point, and the canonical generators of \(U\) by \(\a\) and \(\b\). 
Without loss of generality, we may assume
 that \( s(\s_a - \s_b) = \smax (K) \). From
 Figure~\ref{Fig:Sum}, we see
 that \(\partial ((\s_a - \s_b) \otimes \a) = \s_a \). Since
\(\partial \) is a filtered map of degree \(-1\), we conclude that 
\begin{align*}
s ((\s_a - \s_b) \otimes \a) & \leq s(\s_a) +1 \\
\smax (K) -1 & \leq \smin (K) + 1
\end{align*}
Since we already know that \( \smax (K) \neq \smin (K) \), this gives
the desired result. 
\end{proof}

\subsection{Properties of \(s\)}
We check that \(s\) behaves nicely with respect to mirror image and
connected sum. 

\begin{prop}
\label{Prop:Inverse}
Let \( \overline{K} \) be the mirror image of \(K\). Then we have
\begin{align*}
\smax(\overline{K}) & = - \smin(K)\\
\smin(\overline{K}) & = - \smax(K) \\
s(\overline{K}) & = - s(K)
\end{align*}
\end{prop}
\begin{proof}
Suppose that \(C\) is a filtered complex with filtration \(C=C_0
\supset C_1 \supset \ldots \supset C_n = \{0\} \). Then the dual complex
\(C^*\) has a filtration \(\{0\}= C_0^* \subset C_{-1}^* \subset \ldots
\subset C_{-n}^* = C^* \), where \( C_{-i}^* = \{x \in C^* \ts | \ts
\langle x, y \rangle = 0 , \forall y \in C_i \} \). 

To prove the proposition, we observe that the filtered complex 
\( \ckh ' (\overline{K}) \) is isomorphic to \( (\ckh ' (K))^*
\). Indeed, it is easy 
to see from equation~\ref{Eq:Cobordisms2} that there is an isomorphism
\begin{equation*}
r \co (V, m', \Delta') \to (V^*, \Delta'^*, m'^*) 
\end{equation*}
 which sends
\(\x_\pm \) to \( \x_\mp^* \). Then if \(\x\) is a state of
the diagram \( \overline{K}\), we define \(R(\x) \) to be state of
\(K\) obtained by applying \(r\) all the labels of \(\x\).
It is straightforward to check that the map 
\(R: \ckh ' (\overline{K}) \to (\ckh ' (K))^* \) is the desired
isomorphism. 
(Compare with section 7.3 of \cite{Khovanov}, where it is shown  that 
\( \ckh  (\overline{K}) \cong (\ckh  (K))^* \).)

We now appeal to the following general result, whose
 proof is left to the reader:
\begin{lem}
 If \(C_1\) and \(C_2\)
are dual filtered complexes over a field, then their associated
spectral sequences \(E_n^1 \) and \(E_n^2\) are dual, in the sense that 
\(E_n^1 \cong (E_n^2)^* \). 
\end{lem}
Thus if the two surviving generators in \(E^1_{\infty} \) have
filtration gradings \(\smin \) and \(\smax \), the surviving
generators in \( E^2_\infty \) will have gradings \( - \smax \) and \(
- \smin \). 
\end{proof}

\begin{prop}
\label{Prop:Sum}
\begin{equation*}
s(K_1 \# K_2 ) = s(K_1) + s(K_2)
\end{equation*}
\end{prop}

\begin{proof}
We use the short exact sequence of Lemma~\ref{Lem:SumSeq}. 
Denote the canonical generators of \(K_i\) by \(\s_a^i\) and \( \s_b^i\),
according to their label near the connected sum point. 
It is not difficult to see that \(\kh '(K_1 \# K_2) \) has a canonical
generator \( \s_o\) which maps to \(\s_a \otimes \s_b \) under \(p_*\). 
Thus 
\begin{align*}
s (\s_o) - 1 & \leq s (\s_a^1 \otimes \s_b^2) \\
\smin (K_1 \# K_2) -1 & \leq \smin (K_1) + \smin (K_2)
\end{align*}
Applying the same argument to \(\overline{K}_1\) and \( \overline{K}_2
\), and using the fact that \( \smin (K) = - \smax (K) \), we see that
\begin{align*}
\smax (K_1 \# K_2) + 1 & \geq \smax (K_1) + \smax (K_2) \\
\smin (K_1 \# K_2) + 3 & \geq \smin (K_1) + \smin (K_2) + 4
\end{align*}
Thus
\begin{align*}
\smin( K_1 \# K_2) & = \smin (K_1) + \smin (K_1) + 1 \\
\smax( K_1 \# K_2) & = \smax (K_1) + \smax (K_1) - 1.
\end{align*}
This proves the claim. 
\end{proof}

\section{Behavior under Cobordisms}
\label{Sec:Cobordisms}

Let  \(L_0\) and \(L_1\) be two links in \(\R^3\). An  oriented
 cobordism from  \(L_0\) to
 \(L_1\)  is  a smooth, oriented, compact, 
 properly embedded  surface \(S
\subset \R^3 \times [0,1] \) with \(S \cap (\R^3 \times \{i\}) = L_i \).
 In this section, we define and study a map \( \phi_S \co \kh ' (L_0 )
 \to \kh'(L_1) \) induced by such a cobordism. 
Our construction follows section 6.3 of \cite{Khovanov}, where
Khovanov describes a similar map for the homology theory \( \kh \). 

\subsection{Elementary cobordisms}
Following Khovanov, we  decompose the
cobordism \(S\) into a series of elementary cobordisms, each
represented by a single move from one planar diagram to
another. (See \cite{CarterSaito} for a more detailed treatment of this
material.) For \(i \in [0,1] \), let
\begin{align*}
L_i & = S \cap ( \R^3 \times \{i\}) \\
S_i & = S \cap ( \R^3 \times [0,i]).
\end{align*}
After a small isotopy of \(S\), we can assume that \(L_i\) is a link
in \(\R^3\) for all but finitely many values of \(i\). The orientation
on \(S\) restricts to an orientation on \(S_i\), which in turn
determines an orientation on \(L_i\). We denote this orientation by
\(o_i\). (Note that with this convention, \(o_0\) is the reverse of
the orientation induced on \(L_0\) by \(S\).)

Next, we fix a projection \( p: \R^3 \to \R^2\).
After a further small isotopy of \(S\), we can assume that \(p\)
defines a regular projection of \(L_i\) for all but finitely many
values of \(i\), and that this set of special values is disjoint
from the first set where \(L\) failed to be a link. 
The isotopy type of the oriented planar diagram \(L_i\) remains constant
except when \(L\) passes through one of the two types of 
 special values, where it
changes by some well-defined local move. Each of these moves
corresponds to an elementary cobordism, so we can write the whole
cobordism \(S\) as a composition of elementary cobordisms. 

\begin{figure}
\includegraphics{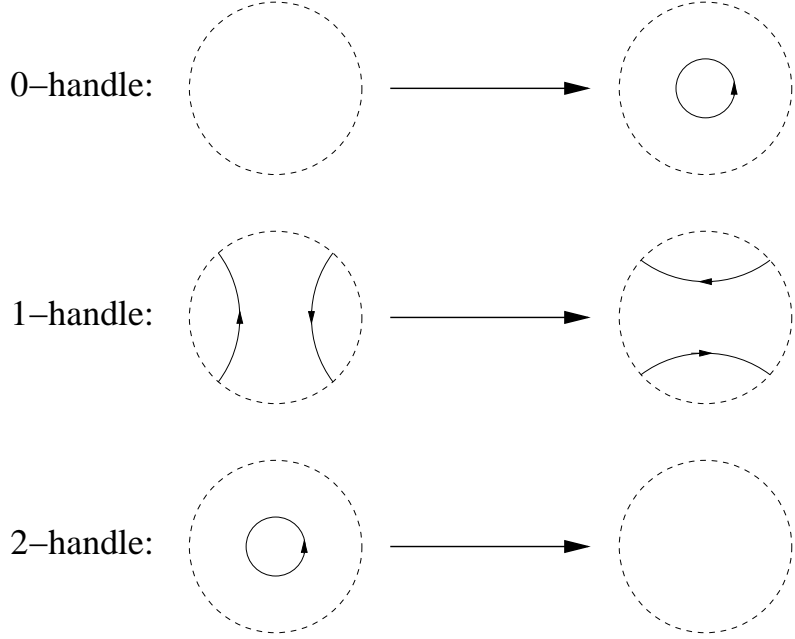}
\caption{\label{Fig:Morse} Local pictures for Morse  moves.}
\end{figure}

The necessary moves may be subdivided into two types:
Reidemeister moves and Morse moves. There is one
Reidemeister-type move for each of the ordinary Reidemeister moves,
 as well as one for each of their inverses. These
moves do not change the topology of the surface \(S_i\).
The Morse moves correspond to the addition of a \(0\), \(1\) or
\(2\)-handle to \(S_i\). They are   illustrated in
 Figure~\ref{Fig:Morse}.

\subsection{Induced Maps}
Given a cobordism \(S\) from \(L_0\) to
\(L_1\), we want to assign to it an induced map \( \phi_S: \kh' (L_0) \to
\kh' (L_1)\) which respects the filtration on \( \kh'\). 
In addition, we would like this assignment to be functorial, in the sense
that if \(S\) is the composition of two cobordisms \(S_1\) and
\(S_2\), \( \phi_S \) is the composition of \( \phi_{S_1} \) and \(
\phi_{S_2} \). Thus it suffices to consider the case when \(S\) is an
elementary cobordism.

Suppose that \( S\)  is an elementary cobordism corresponding 
 to the \(i\)-th Reidemeister move or its
 inverse. Then  we define  \( \phi_S \) to be \( \rho_{i*}'\) or its
 inverse. By Proposition~\ref{Prop:ReidFilt}, this is a filtered map of
 degree \(0\). 
If \(S\) is an elementary cobordism corresponding to a
Morse move, then we take \( \phi_S \)  to be the map  induced
 by \(\psi : \ckh' (L_0) \to \ckh' (L_1) \), where \( \psi \) is the 
 result of applying the 
TQFT \( \A '\) to the corresponding map of cubes.
 In other words, if the move corresponds to the addition of a
 \(0\)-handle or a \(2\)-handle, we apply
\(\iota ' \) or \( \epsilon '\), respectively, to the summand at each
 vertex of the cube. If it corresponds to the addition of a
 \(1\)-handle, we apply either \( m'\) or \( \Delta '\), depending
 on whether the move results in a merge or a split at the vertex in
 question.
 It is easy to see that \( \phi_S\) is a
 filtered map of degree \(1\) for a \(0\)-- or
  \(2\)--handle addition and degree \(-1\) for a \(1\)--handle. 

In general, given a cobordism \(S\), we decompose it as a union of 
elementary cobordisms: \(S=S_1 \cup S_2  \ldots
\cup S_k \) and
  define the induced morphism \(\phi_S : \kh' (L_0) \to \kh ' (L_1) \)
to be the composition \( \phi_{S_k} \circ \ldots \circ \phi_{S_1} \),
which is a 
filtered map of degree \( \chi (S) \). 
We expect that the map \(\phi_S\) will depend only on the isotopy
class of \(S \ \text{rel} \ \partial S\) 
({\it c.f} \cite{Magnus}, where an analogous result
is proved for the Khovanov homology), 
but since we do not need this fact, we will not pursue it here. 

\subsection{Canonical generators}

The maps \( \phi_S\)
 behave nicely with respect to  canonical generators.
\begin{prop}
\label{Prop:Maps}
Suppose \(S\) is an oriented  cobordism from \(L_0\) to \(L_1\) which
is \emph{weakly connected}, in the sense that every component of \(S\)
has a boundary component in \(L_0\). Then \( \phi_S([\s_{o_0}])\) is a
nonzero multiple of \([\s_{o_1}] \).
\end{prop}

\noindent {\bf Remark:} Some sort of  connectedness hypothesis is clearly
necessary for the 
proposition to hold. For example, if we take \(S\) to be the 
union of a product cobordism and a trivially embedded sphere, the
induced map on \(\kh'\) is  the zero map. 

\begin{proof}
In fact, we will prove a slightly stronger statement. Suppose \(i\) is
a regular value for the cobordism \(S\), so that \(L_i\) is a link. We
divide the components of \(S_i\) into two sorts: those of the {\it
  first type}, which have a boundary component in \(L_0\), and those
of the {\it second type}, which do not. We say that an orientation \(o\)
on \(S_i\) is {\it permissible} if it agrees with 
the  orientation of \(S\) on components of the first type.
(Here and in what follows, we use \(o_I\) to denote both a permissible
 orientation on \(S_i\) and the orientation it induces on \(L_i\).)
 We claim that 
\begin{equation*}
\phi_{S_i}([\s_{o_0}]) = \sum_{I} a_I [\s_{o_I}]
\end{equation*}
where \(\{o_I\}\) runs over the set of permissible orientations on \(S_i\)
and each coefficient \(a_I\) is nonzero.
Note that the weak connectivity hypothesis 
 implies that there is only one permissible
orientation on \(S_1\), so the proposition is implied by the claim. 

To prove the claim, it suffices to check that if it holds for
\(S_i\), then it holds for \(S_{i'}\) as well, where \(S_i'\) is
the composition of \(S_i\) with a single elementary cobordism
\(S_e\). If this cobordism
corresponds to a Reidemeister type move, this is a straightforward
consequence of Proposition~\ref{Prop:ReidCan}. Below, we check that it
holds for each of the Morse-type moves as well.
\vskip0.05in
\noindent{\it \(0\)-Handle Move:} In this case, \(\phi _{S_{e}} (\s_{o_I}) =
  \s_{o_I} \otimes \frac{1}{2}(\a-\b) \), where the second factor in
  the tensor product refers to the labels on the newly created circle. 
  \(S_{i'}\) has a new component of the second type --- namely, the
  disk bounded by the new circle --- and \(\s_{o_I} \otimes \a \) and
  \(s_{o_I} \otimes \b \) are the canonical generators corresponding
  to the two possible orientations on \(S_{i'} \) which agree with 
\(o_I\) on all components other than the new one. 

\vskip0.05in
\noindent{\it \(1\)-Handle Move:} Suppose that the orientation \(o_I\) is
  actually the orientation \(o_i\)
 induced by \(S_i\). Then  the two strands involved in the move
  have opposite orientations, so by Lemma~\ref{Lem:CO}, they must have
  the same label. Since
\begin{align*}
m'(\a\otimes \a) & = 2 \a  \qquad  &\Delta ' (\a) & = \a \otimes \a \\
m'(\b\otimes \b) & = -2 \b  \qquad  &\Delta ' (\b) & = \b \otimes \b
\end{align*}
we see that
\( \phi_{S_e}(\s_{o_i}) \) is a nonzero multiple of \( \s_{o_{i'}} \),

More generally, the orientation \(o_I\) is either compatible with some
orientation \(o_e\) on \(S_e\), or it is not. In the former case, the two
strands involved in the move point in opposite directions and have the
same label, and \( \phi_{S_e}(\s_{o_I}) \) is a nonzero multiple of 
\( \s_{o_I'} \) where \(o_I'\) is the orientation induced on
\(L_{i'}\) by \(o_e\). In the latter case, the two strands point in the same
direction and have different labels, so \( \phi_{S_e}(\s_{o_I}) =
0\). 

Now we consider what happens to the components
of \(S_i\) during the move. If the move splits one component of
\(L_i\) into two components of \(L_{i'} \), then the number and type
of components of \(S_i\) remains constant. In this case, the set of
permissible orientations on \(S_i\) is naturally identified with the set
of permissible orientations on \(S_{i'}\).  There is always an
orientation on \(S_{e} \) compatible with \(o_I\), and \( \phi_{S_e}
(\s_{o_I}) \) is a nonzero multiple of \(\s_{o_I'} \). 

On the other hand, if the move merges two components of 
\(L_{i}\) into one component of \(L_{i'} \), there are several
possibilities to consider. If the merge involves only a single component
of \(S_i\), the situation is like  the one above: there is always an
orientation on \(S_{e} \) compatible with \(o_I\), and \( \phi_{S_e}
(\s_{o_I}) \) is a nonzero multiple of \(\s_{o_I'} \). The same
argument applies when \(S_e\) merges two components of \(S_i\), both
of which are of the first type. 

Finally,
suppose the merge joins two components of \(S_i\), at least one of
which is of the second type. Then the set of permissible orientations on
\(S_{i'}\) is only half as large as the set of permissible
orientations on \(S_i\). If \(o_I\) extends to a permissible
 orientation \(o_I'\) on  \(S_{i'} \), \( \phi_{S_e} (\s_{o_I}) =
 \s_{o_I'} \), while if it does not, \( \phi_{S_e}(\s_{o_I}) = 0 \). 

\vskip0.05in
\noindent{\it \(2\)-Handle Move:} In this case, a permissible orientation \(o_I
\) on \(S_i\) extends to a unique permissible orientation \(o_I'\) on
\(S_{i'}\). Since \(\epsilon'(\a) = \epsilon'(\b) = 1\), \( \phi_{S_e}
(\s_{o_I}) = \s_{o_I'}\). To prove the claim, it suffices to show that
two permissible orientations on \(S_i'\) cannot induce the same
orientation on \(L_{i'} \). But if this were the case, \(S_i\) would
 have a closed component, contradicting the hypothesis that
\(S\) is weakly connected. 

\end{proof}

\begin{cor}
If \(S\) is a connected cobordism between knots \(K_0\) and \(K_1\),
 then \( \phi_S \) is an isomorphism. 
\end{cor}
\begin{proof}
Fix an orientation \(o\) on \(S\). Then \(\{\s_{o_0},
\s_{\overline{o}_0} \} \) is a basis for \( \kh' (K_1) \). Its image
under \( \phi_S \) is \(\{k_1\s_{o_1},
k_2 \s_{\overline{o}_1} \} \) \((k_1, k_2 \neq 0 ) \),
 which is a basis for \( \kh ' (K_2) \). 
\end{proof}

\subsection{The slice genus}

 We can now prove the first two theorems from the introduction.

\begin{proof} (of Theorem~\ref{Thm:One}.)
Suppose \(K \subset S^3 \) bounds an oriented surface of genus \(g\) in \(
B^4\). Then  there is an orientable connected cobordism of 
Euler characteristic \(-2g\)  between 
\(K   \) and the unknot \(U\)  in \( \R^3 \times [0,1] \). Let \( x \in
\kh' (K) -\{0\} \) be a class for which \( s (x) \) is maximal. Then
\( \phi_S(x) \) is a nonzero element of \( \kh ' (U) \). Now
\( \phi_S\) is a filtered map with filtered degree \(-2g\), so 
\begin{equation*}
 s(\phi_S(x)) \geq s(x) - 2g. 
\end{equation*}
 On the other hand, \( \smax (U) =
1\), so
\begin{equation*}
  s( \phi_S(x)) \leq 1. 
\end{equation*}
It follows that \( s(x) \leq 2g+1
\), so \( \smax (K) \leq 2g+1\) and \( s(K) \leq 2g \). To show that
\( s(K) \geq -2g \), we 
apply the same argument to \(\overline{K} \)  (which bounds a
surface \( \overline{S} \) of genus \(g\)) and use the fact that 
\( s (\overline{K}) = - s (K)\). 
\end{proof}

\begin{proof}(of Theorem~\ref{Thm:Two}.)
If \(K_1\) and \(K_2\) are concordant, then \(K_1 \# \overline{K_2}\) is slice,
so
\begin{equation*}
0 = s(K_1 \# \overline{K_2}) = s(K_1) - s(K_2).
\end{equation*}
Thus \(s\) gives a well-defined map from 
\( \text{Conc} (S^3)\) to  \( \Z \). 
That this map is a homomorphism is immediate from
Propositions~\ref{Prop:Inverse} and \ref{Prop:Sum}. 
\end{proof}

\begin{cor}
Suppose \(K_+\) and \(K_-\) are knots that differ by a single crossing
change --- from a positive crossing in \(K_+\) to a negative one in \(
K_-\). Then
\begin{equation*}
 s(K_-) \leq s(K_+) \leq s(K_-) + 1
\end{equation*}
\end{cor}

\begin{proof}
In \cite{Livingston}, Livingston shows that this skein inequality
holds for any knot invariant satisfying the properties of
Theorems~\ref{Thm:One} and \ref{Thm:Two}. 
\end{proof}

\section{Computations and Relations with other Invariants}

Although the invariant \(s(K)\) is algorithmically computable from a
diagram of \(K\), it is impossible to compute by hand for all but the
smallest knots. In this section, we describe some techniques which enable
us to efficiently compute \(s\).

\subsection{Using \(\kh\)}

For many knots, it is a  simple matter to compute \( s(K)
\) from the ordinary Khovanov homology \( \kh (K) \). Although \( \kh
(K) \) is also hard to compute by hand, there are already
 a number of computer programs available for
this purpose, including Bar-Natan's pioneering program \cite{DBN}
 and a more
recent, faster program written by Shumakovitch \cite{Shumak}.

In \cite{DBN}, Bar-Natan made the following observation,
 based on his computations of \(\kh \) for knots with 
\(10\) and fewer crossings. 

\begin{conj} (Bar-Natan) 
The graded Poincare polynomial \( P_{\kh}(K) \) of \( \kh (K) \) has
the form
\begin{equation*}
P_{Kh} (K) = q^{s(K)}(q+ q^{-1}) + (1 + tq^4) Q_{Kh}(K)
\end{equation*}
where \(Q_{Kh} (K) \) is a polynomial with all positive coefficients. 
\end{conj}

\noindent In \cite{ESL2}, Lee showed that this conjecture holds whenever
her spectral sequence for \( \kh ' \) converges after the \(E_2\)
term. In this case, it is easy to see that the invariant \(s(K)\) is
equal to the exponent \(s(K) \) which appears in Bar-Natan's conjecture.

 To see how widely applicable this condition is, we introduce
the notion of the homological {\it width}  of a knot.

\begin{dfn}
If \(K\) is a knot, let \( \mu (K) = \{a - 2b \ts | \ts q^a t^b \text
  {is a monomial in } P_{Kh}(K) \} \). The  width  \(W(K)\) is one more
  than the
  difference between the maximum and minimum elements of \( \mu (K)
  \). 
\end{dfn}

\noindent In other words, \(W (K) \) is the number of diagonals  in the convex
hull of the support of \(\kh (K) \). 

\begin{prop}
If \(W(K) \leq 3\), then the spectral sequence for \(\kh' (K) \)
converges after the \(E_2\) term, and our \(s(K) \) is the same as
Bar-Natan's. 
\end{prop}

\begin{proof}
Suppose \(W(K) \) has width \( \leq 3\). Then if \(x\) is an element
of \( \kh ' (K) \) with \(q\)-grading \(a\) and homological grading
\(b\), the minimum possible \(q\)-grading of an element with
homological grading \(b-1\) is \(a-6\). Since the differential \(d_n\)
 on the \(E_n\) term of the spectral sequence lowers the \(q\)-grading by
\(4(n-1)\), \(d_n\) must be trivial for all \(n \geq 3\). 
\end{proof}

\noindent Theorem~\ref{Thm:Three} follows from this fact, since 
 Lee has shown  \cite{Lee} that if \(K\) is an alternating knot, then
it has width two and Bar-Natan's
\(s\) is equal to \( \sigma (K) \).  

The proposition also applies to many non-alternating knots. 
Indeed, 
using Shumakovitch's tables and a computer, it is
straightforward to check that there are only four knots with
13  or fewer crossings whose width is greater than three. 
Inspecting \(\kh\) of these four exceptions, one sees that in each
case, the spectral sequence must converge after the \(E_2\) term. 
 Thus for all knots with 13 or fewer crossings, the value of \(s(K) \) agrees
with the value of Bar-Natan's \(s\)  tabulated in \cite{DBN}
and \cite{Shumak}. Below, we list those knots of 11 crossings or fewer
for which \(s(K) \neq \sigma (K) \). There are 22 such knots, and
\(|s(K) | > |\sigma (K) | \) (and thus provides a better bound on the
slice genus) for precisely half of them.

$$
\begin{array}{|l|r|r||l|r|r||l|r|r|}
\hline
\rule{0pt}{14pt} K & s(K) & \sigma(K) &
 K & s(K) & \sigma(K) &
 K & s(K) & \sigma(K)  \\
\hline
\hline
9_{42} & 0 & 2 &
11_{n9} & 6 & 4 &
11_{n70} & 2 & 4 
\\
10_{132} & -2 & 0 & 
11_{n12} & 2 & 0 &
11_{n77} & 8 & 6 
\\
10_{136} & 0 & 2 &
11_{n19} & -2 & -4 &
11_{n79} & 0 & 2 
\\
10_{139} & 8 & 6 &
11_{n20} & 0 & -2 &
11_{n92} & 0 & -2 
\\
10_{145} & -4 & -2 &
11_{n24} & 0 & 2 &
11_{n96} & 0 & 2 
\\
10_{152} & -8 & -6 & 
11_{n31} & 4 & 2 &
11_{n138} & 0 & 2 
\\
10_{154} & 6 & 4 &
11_{n38} & 0 & 2 &
11_{n183} &  6 & 4
\\
10_{161} & -6 & -4 &
& & &
& & \\
\hline
\end{array}
$$
\vskip0.07in
Knots with 10 or fewer crossings are labeled according to their
numbering in Rolfsen, while those with 11 crossings use the {\it
  Knotscape} numbering. The values of the signature are taken from
\cite{DBNAtlas}. 
All of the knots in the table have a homological width of 3,
which raises the following question:
if \(K\) has homological width 2 ({\it i.e.} is H-thin
in the terminology of \cite{Khovanov2}), must \(s(K) = \sigma (K)\)? 
\vskip0.03in

\subsection{Positive knots}
If \(K\) is a positive knot, \(s(K) \) can be computed directly from
the definition. To see this, consider a canonical generator
\(\s_o\) for a positive diagram of \(K\). Since each crossing of
 \(K\) is positive, its oriented resolution is the
\(0\)-resolution. Thus the state
\(\s_o\) lives in the extreme corner of the cube of resolutions: it
has homological grading \(0\), and there are no generators in \( \ckh'
(K) \) with homological grading \(-1\). It follows that the only class
homologous to \( \s_o\) is \(\s_o\) itself, so
\begin{equation*}
\smin (K) = s([ \s_o]) = q(\s_o)
\end{equation*}

To compute \(q(\s_o) \), we change back to the basis \(\{ \x_-, \x_+\}
\). In the expansion of \(\s_o\) with respect to this basis, there is a
unique state with minimal \(q\)-grading, namely, 
the state in which every circle of the oriented
resolution is labeled with a \( \x_-\). If the positive diagram of \(K\)
has \(n\) crossings, and its oriented resolution has \(k\) circles, then
\begin{align*}
 q ( \s_o) & = p(\s_o) + \gr (\s_o) + n_+ - n_- \\
& = -k + 0  + n -0 
\end{align*}
so
\begin{equation*}
s (K) = -k + n + 1
\end{equation*} 

On the other hand, Seifert's algorithm gives a Seifert surface \(S\) for
\(K\) with euler characteristic \(k - n\), so 
\begin{equation*}
2g(K)  \leq 2 g(S) = n-k + 1 = s(K) \leq 2 g_*(K)
\end{equation*}
Since \(g_*(K) \leq g(K) \), the inequalities above must all be
equalities. This completes the proof of Theorem~\ref{Thm:Four}. 

\subsection{Comparison with \(\tau\)}
We end this section by commenting on the conjecture relating \(s\) and
\(\tau\) which was stated in the
introduction.  In addition to the fact that the two invariants share
 the  properties of Theorems 1
through 4, there is a good deal of
numerical evidence supporting the conjecture. Recently, a fair amount of work has been done on the problem of
computing \( \tau \) for knots with 10 and fewer crossings. Combining
the results of  \cite{Goda},
 \cite{Livingston}, \cite{OS8}, \cite{OS7}, and \cite{OS10}
 with some unpublished computations of the
author, it appears that the value of \( \tau \) has been determined
for all but two 
knots of 10 crossings and fewer. (The exceptions are
 \(10_{141} \) and \(10_{150}\).)
 For all of these knots, \( s = 2 \tau \). The equality can
also be checked on certain special classes of knots, such as the
pretzel knots of \cite{OSMut}. 
 If the conjecture were true, it would make many computations in knot
 Floer homology  easier. (For example, with our current technology, 
 it seems like quite a
 laborious project to compute \( \tau \) for all 11-crossing
 non-alternating knots.) Even if it is not true, we hope that
 the remarkable similarity between the two theories will have an enlightening
 explanation.

\section{Reidemeister Moves}
\label{Sec:Reidemeister}
In this section, we prove  the results involving
Reidemeister moves which were stated in section~\ref{Sec:Review} and
\ref{Sec:Def}. 

\begin{proof} (of Theorem~\ref{Thm:Lee1}.)
The proof that the desired  spectral sequence exists was sketched in 
section~\ref{Sec:Review}. 
To prove its invariance, we use the following
basic lemma, whose proof may be found in \cite{UsersGuide}, Proposition 3.2.
\begin{lem}
\label{Lem:SSInv}
Suppose \(F \colon C_1 \to C_2 \) is a map of filtered complexes which
respects the filtrations. Then \(F\) induces maps of spectral
sequences \(F_n \colon E^1_n \to E^2_n \), and if \(F_n \) is an
isomorphism, \(F_m\) is an isomorphism for all \(m \geq n\). 
\end{lem}

In section 4 of \cite{ESL2}, Lee proves the invariance of \( \kh'\) by checking
its invariance under the three Reidemeister moves. For each move, she
exhibits a chain map between the complexes associated to the link
diagram 
before and after the move. To prove the theorem, it suffices to
check that these maps respect the \(q\)-filtration, and that they
induce isomorphisms on the \(E_2\) terms. The latter
claim is straightforward, since in each case the induced maps on the
\(E_1\) terms are identical to the maps used in section 5 of
\cite{Khovanov} to prove invariance of \(\kh\). Below, we sketch
the proof of invariance for each move and explain why the maps in
question respect the filtrations. For full details, we refer the
reader to \cite{Khovanov} and \cite{ESL2}.
\vskip0.05in
\begin{figure}
\includegraphics{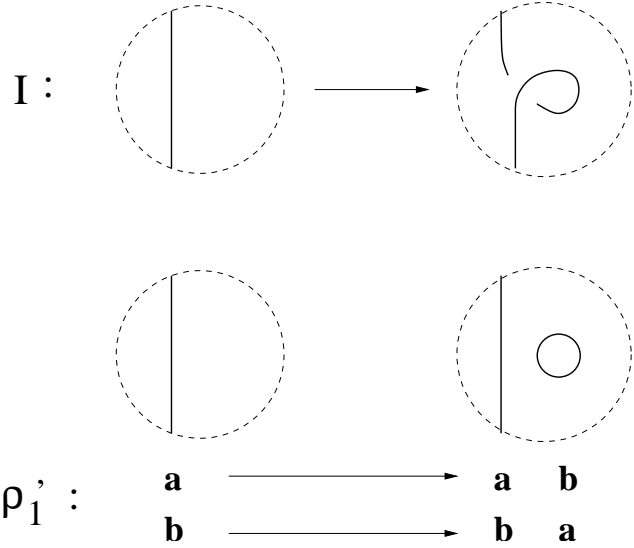}
\caption{\label{Fig:Reidemeister1} The Reidemeister I move and 
the map \( \rho'_1\).}
\end{figure}
\noindent{\it Reidemeister I Move:} Let \(\tilde{L}\) be the diagram
\(L\) with an additional left-hand curl added in. Then
\(\ckh'(\tilde{L})\) can be decomposed as a direct sum \(X_1 \oplus
X_2\), where \(X_2\) is acyclic and \(X_1\) is isomorphic to \( \ckh '
(L) \) via the map \(\rho_1'\co \ckh' (L) \to X_1 \) illustrated in
Figure~\ref{Fig:Reidemeister1}. 
 In terms of the basis \(\{\x_\pm\}\), we have
\begin{align*}
\rho_1 ' (\x_-) & = \x_- \otimes \x_- - \x_+ \otimes \x_+ \\
\rho_1 ' (\x_+) & = \x_+ \otimes \x_- - \x_- \otimes \x_+ 
\end{align*}
The corresponding map \( \rho_1\) in \cite{Khovanov} is given by
\begin{align*}
\rho_1  (\x_-) & = \x_- \otimes \x_-  \\
\rho_1  (\x_+) & = \x_+ \otimes \x_- - \x_- \otimes \x_+ 
\end{align*}
so \(\rho_1 '\)  is filtration non-decreasing, and its induced map on
\(E_1\) terms is \( \rho_1\).

\noindent {\bf Remark:} There is another version of the first
Reidemeister move, corresponding to the addition of a right-hand curl.
Although it is not difficult to define an appropriate map \( \rho_{1'}'\) for
  this move directly, for the sake of brevity we adopt the solution of
  \cite{DBN} and \cite{ESL2} and  define it to be the
  composition of maps induced by an appropriate Reidemeister II move
  followed by a Reidemeister I move. 
\vskip0.05in
\noindent{\it Reidemeister II Move:} Let \(L\) and \(\tilde{L}\) be as
shown in figure~\ref{Fig:Reid2}. In this case, 
\(\ckh'(\tilde{L})\) can be decomposed as a direct sum \(X_1 \oplus
X_2\oplus X_3\), 
where \(X_2\) and \(X_3\)  are acyclic and there is an isomorphism
\( \rho_2 ' : \ckh' (L) \to X_1\), which is given by 
\begin{equation*}
\rho_2 ' (z) = (-1)^{\gr(z)}(z + \iota (d'_{01\to 11} (z)))
\end{equation*}
The maps \(\iota \) and \(d'_{01\to 11} \) are shown in the
figure. 
The isomorphism \( \rho_2 \) in \cite{Khovanov} has the same form, but with
\(d_{01\to 11} \) in place of \( d'_{01\to 11}\). Since 
\(d-d'\) is strictly filtration increasing, it follows that 
\(\rho_2 ' \) is filtration non-decreasing, and its induced map on
\(E_1\) terms is \( \rho_2\). 

\begin{figure}
\includegraphics{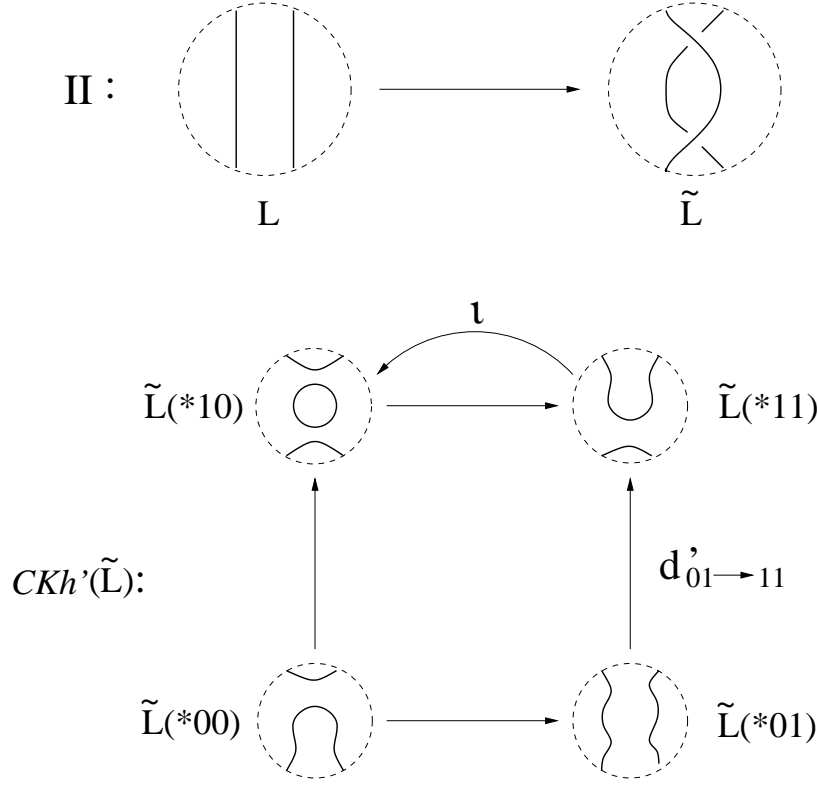}
\caption{\label{Fig:Reid2} The Reidemeister II move and the maps \(\iota \) and \( d'_{01\to 11}
\).}
\end{figure}
\vskip0.05in
\noindent{\it Reidemeister III Move:}
\begin{figure}
\includegraphics{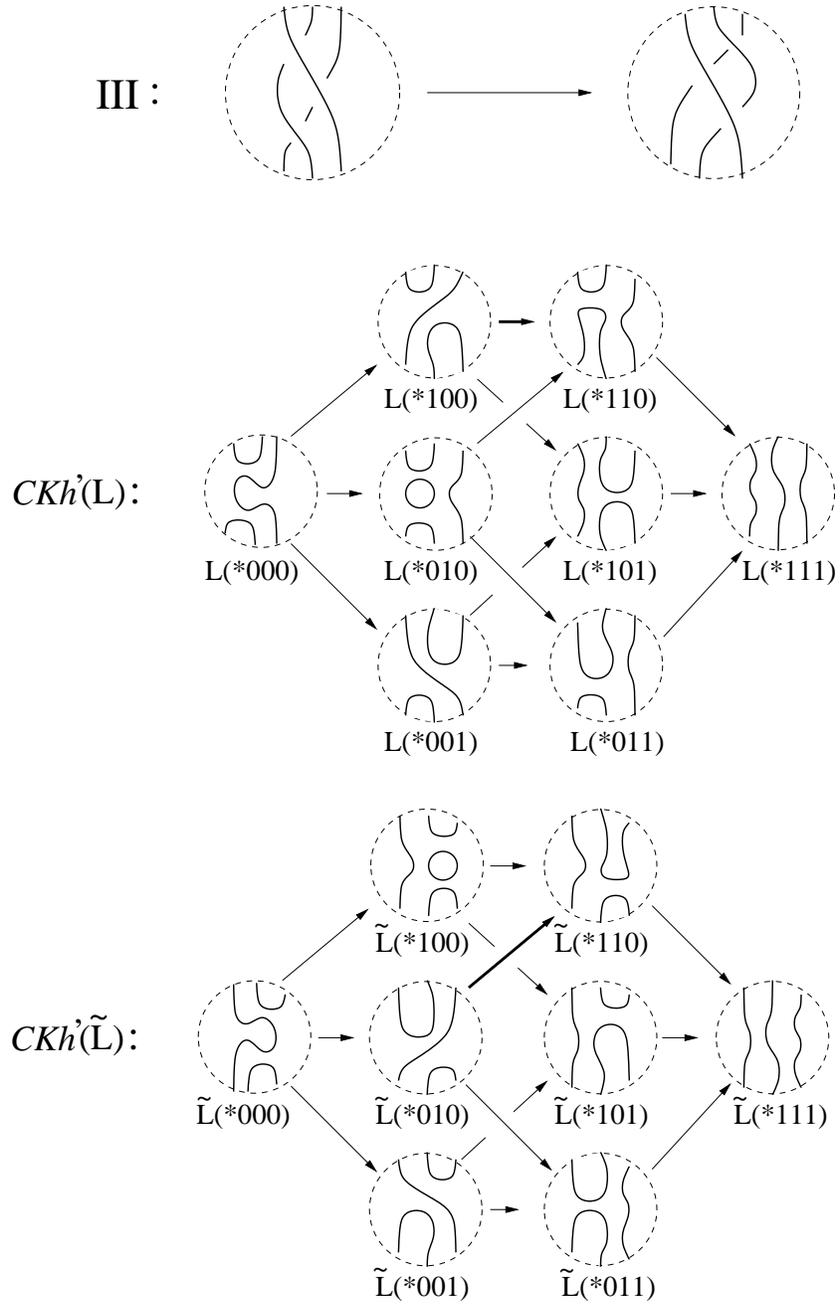}
\caption{\label{Fig:Reid3} The Reidemeister III move. The relevant
  components of the differentials (\(d'_{100\to110} \) and 
\(d'_{010\to110}\)) are marked in bold.  }
\end{figure}
Let \(L\) and \( \tilde L\) be as shown in Figure~\ref{Fig:Reid3}.
Then there are direct sum decompositions
\begin{align*}
\ckh'(L) & \cong X_1 \oplus X_2 \oplus X_3 \\
\ckh'(\tilde{L}) & \cong \tilde{X}_1 \oplus \tilde{X}_2 \oplus \tilde{X}_3
\end{align*}
where \(X_2,X_3,\tilde{X}_2\), and \(\tilde{X}_3\) are acyclic and
there is an isomorphism \( \rho_3 '\co X_1 \to \tilde{X}_1 \). To describe
 \(X_1\) and \( \tilde{X}_1 \), we first define maps
\begin{align*}
\beta' \co & \ckh '(L(*100)) \to \ckh'(L(*010)) \\
\tilde{\beta}'\co & \ckh '(\tilde{L}(*010)) \to \ckh'(\tilde{L}(*100))
\end{align*}
by 
\begin{align*}
\beta' & = \iota \circ d'_{100\to110} \\
\tilde{\beta}' & = \iota \circ d'_{010\to110}
\end{align*}
Then
\begin{align*}
X_1  & = \{x + \beta'(x) + y \ts | \ts x \in \ckh '(L(*100)), y\in \ckh
'(L(*1))\} \\
\tilde{X}_1  & = \{x + \tilde{\beta}'(x) + y \ts | \ts x \in \ckh '
(\tilde{L}(*010)), y\in \ckh'(\tilde{L}(*1))\}
\end{align*}
and 
\begin{equation*}
\rho_3'(x+\beta'(x)+y) = x + \tilde{\beta}' (x) + y. 
\end{equation*}
The isomorphism \( \rho_3 \) in 
\cite{Khovanov} is defined similarly, except that it uses \(d\) instead of
\(d'\) to define maps  \( \beta \) and \( \beta '\). Since \(d'\) does
not increase the \(q\)-grading, we clearly have \(q(\beta'(x)) \geq 
q (x) \). From this, it follows that \(\rho_3 '\) does not decrease the
\(q\)-grading. Since \(d-d'\) strictly increases the \(q\)-grading,
 the map induced on \(E_1\) terms by \( \rho_3 '\) is equal to \( \rho_3
 \). To finish the proof, we apply  Lemma~\ref{Lem:SSInv} three times:
 first to the inclusions \(X_1 \into \ckh '(L) \) and
 \(\tilde{X}_1 \into \ckh '(\tilde{L}) \), and then to the map  \( \rho_3 '\). 

\end{proof}

\begin{proof} (of Proposition~\ref{Prop:ReidCan}.)
We check the claim directly for each  Reidemeister move:
\vskip0.05in
\noindent{\it Reidemeister I Move:}
In this case, it is easy to see that \( \rho _1' (\s_o) =
\s_{\tilde{o}} \). 
\vskip0.05in
\noindent{\it Reidemeister II Move:}
Suppose that the two strands in \(L\) point in the same
direction. Then by Lemma~\ref{Lem:CO}, they have different labels, so 
\(d'_{01\to11} (\s_o) = 0\). The oriented resolution of \(\tilde{L} \)
is contained in  \(\ckh'(\tilde{L}(*01)) \cong
\ckh' (L)  \), so \(\rho_2' (\s_o) = (-1)^0(\s_{\tilde{o}}) =
\s_{\tilde{o}}\).

Now suppose the two strands point in different directions, so that they
have the same label. Let us assume for the moment that this label is \(\a\).
Then we define \(\s_{\tilde{ij}} \in \kh' (\widetilde{L}(*ij)) \) be the
state which is identical to \(\s_o\) outside the area where the move
takes place and has all components inside the area of the move labeled
with an \(\a\). Then a direct computation shows that either
 \begin{align*}
\rho_2'(\s_o) & = \s_{\widetilde{01}} + \frac{1}{2}(\s_{\widetilde{10}}
- \s_{\tilde{o}}) \\ 
&  =  - \frac{1}{2} ( \s_{\tilde{o}} +  d' (\s_{\widetilde{00}}))
\end{align*}
 if the two strands belong to
the same component, or
 \begin{align*}
\rho_2'(\s_o) & = \s_{\widetilde{01}} + (\s_{\widetilde{10}}
- \s_{\tilde{o}}) \\ 
&  =  -  ( \s_{\tilde{o}} +  d' (\s_{\widetilde{00}}))
\end{align*}
 if they belong to different components. This proves the claim in the
 case where both strands  are labeled with an \(\a\). We leave it to
 the reader to check that a similar argument applies when they are
 both labeled with a \(\b\). 

\vskip0.05in
\noindent{\it Reidemeister III Move:}
\begin{figure}
\includegraphics{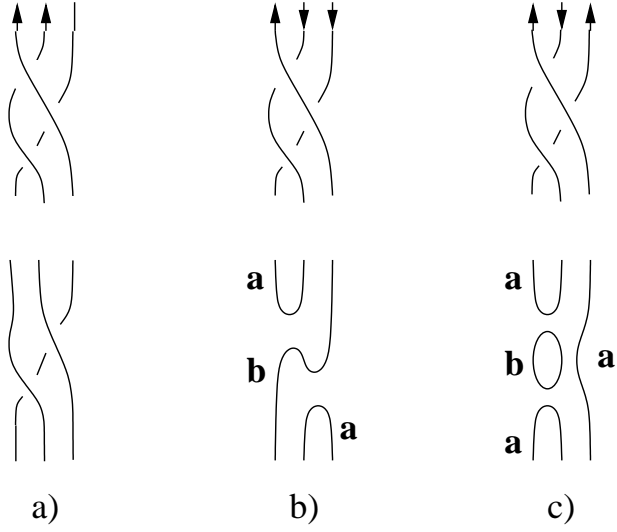}
\caption{\label{Fig:Reid3a} Possible orientations for \(L\) and their
  respective canonical generators.}
\end{figure}
Here there are three cases to consider. First, suppose that
 the two overlying strands in \( L\)
are oriented as shown in Figure~\ref{Fig:Reid3a}a. Then \(\s_o \in \ckh
' (L(*1)) \), and it is easy to see that \( \rho'_3 (\s_o) =
 \s_{\tilde{o}} \).

Next, suppose that the three strands are  oriented as shown in 
Figure~\ref{Fig:Reid3a}b. Then \(\s_o \in 
\ckh'(L(*100)) \) and \( \s_{\tilde{o}} \in \ckh' (\tilde{L}(*010)) \).
Clearly \( \beta ' (\s_o) = \tilde{\beta} ' (\s_{\tilde{o}}) = 0 \), 
  so \(\s_o \in X_1\) and \( \s_{\tilde{o}} \in \tilde{X}_1 \). Again,
  it follows that \(\rho_3'(\s_o) = \s_{\tilde{o}} \). 

Finally, suppose the strands are oriented as shown in
Figure~\ref{Fig:Reid3a}c. In this case, the oriented resolution of
\(L\) is in \(L(*010)\), and the oriented resolution of 
\(\tilde{L}\) is in \( \tilde{L}(*100) \). Inside the region under
consideration, \(\s_o\) looks like the state of
Figure~\ref{Fig:Reid3a}c
 (perhaps with \( \a\)'s and \(\b\)'s reversed.)
Our first step is to exhibit some \( \mathfrak{t} \in X_1\) which
is homologous to \(\s_o\).
As before, we denote by \(\s_{ijk} \) the unique state of \(L(*ijk) \)
which is the same as \(\s_o\) outside the area of the Reidemeister
move and has all its components inside this area labeled by \(\a\)'s.

Assume for the moment that
all three strands  shown in \(L(*000) \) belong to different
components. In this case, we can take 
\begin{equation*}
\mathfrak{t} = \s_o - 2 \s_{100} - \s_{010} - 2 \s_{001} = \s_o - d'
(\s_{000}). 
\end{equation*}
Indeed,  \( \beta '(- 2 \s_{100}) = \s_o -
\s_{010} \) and \(  \s_{001}  \in \ckh ' (L(*1)) \), so \(
\mathfrak{t} \in X_1\). Then 
\begin{align*}
\rho'_3(\mathfrak{t}) & = -2 \s_{\widetilde{010}} - 2
\tilde{\beta}'(\s_{\widetilde{010}}) - 2\s_{\widetilde{001}} \\
 & = -2\s_{\widetilde{010}} - 2 \s_{\widetilde{100}} + 2 \s_{\tilde{o}} -
2\s_{\widetilde{001}} \\
& = 2 \s_{\tilde{o}} - d' (\s_{\widetilde{000}})
\end{align*}
which proves the claim. 

We leave it to the reader to check that a similar argument applies 
to each of the four other ways in which the segments outside the area of
the move can be connected, as well as  when the roles of \(\a\)
and \(\b\) are reversed. In each case, it is not difficult to verify
that \(\rho'_{3*}([\s_0]) \) is one of \( \pm [\s_{\tilde{o}}], \pm 2
[\s_{\tilde{o}}]\), or \( \pm \frac{1}{2} [\s_{\tilde{o}}]\).

\end{proof}

\begin{proof} (of Proposition~\ref{Prop:ReidFilt}.)
In the case of \( \rho_{1*}' \) and \( \rho_{2*}'\), the claim is
immediate, since these maps are induced by filtered chain maps. For
the others, we use the following

\begin{lem}
Suppose \(f\co C_1 \to C_2 \) is a map of filtered chain complexes
with the property that the induced map of spectral sequences 
 \(f_2 \co E_1^2 \to E_2^2 \) is an
 isomorphism. Then \(f_*^{-1}\) is a filtered map with respect to the
 induced filtrations on \(H_*(C_1)\) and \(H_*(C_2) \).
\end{lem}
\begin{proof}
Since \(f_2\) is an isomorphism, \(f_{\infty}\) (the induced map on
filtered gradeds) is as well. It follows that \(f_*\) is an
isomorphism. Suppose \(f_*^{-1}\) does not respect the
filtration. Then there must be some \( \x \in H_*(C_1) \) whose
filtration is strictly increased by \(f_*\). But this contradicts the
fact that \( f_{\infty} \) is an isomorphism. 
\end{proof}
The remaining cases now follow easily from the results used in the
proof of Theorem~\ref{Thm:Lee1}. Indeed, \( \rho_1'\) and \(\rho_2'\)
both induce isomorphisms of \(E_2\) terms, and \( \rho_{3*}' = \iota_{1*}
\circ \psi_* \circ \iota_{2*}^{-1} \), where \( \iota _1\), \(\iota_2\), and \(
\psi \) all induce isomorphisms of \(E_2\) terms. 
\end{proof}

\bibliographystyle{plain}
\bibliography{../mybib}

\end{document}